\documentclass[12pt]{article}
\input epsf
\usepackage{epsfig}
\usepackage{amsmath,amsfonts,amsthm,bm}
\begin{document}
\def\C{\mathbf{C}}
\def\bC{\mathbf{\overline{C}}}
\def\R{\mathbf{R}}
\def\Sym{\mathrm{Sym}\, }
\def\Rea{\mathrm{Re}\, }
\def\RP{\mathbf{RP}}
\title{Co-axial monodromy}
\author{Alexandre Eremenko\thanks{Supported by NSF grant DMS-1665115.}}
\maketitle
\begin{abstract} For Riemannian metrics of constant positive curvature
on a punctured sphere with conic singularities at the punctures
and co-axial monodromy of the developing map, possible
angles at the singularities are completely described.
This completes the recent result
of Mondello and Panov.

The related problem of describing
possible multiplicities of critical points of logarithmic potentials
of finitely many charges is also solved. 
\vspace{.1in}

\noindent
2010 MSC: 57M50, 53C45, 31A99. Keywords: conic singularities,
positive curvature, developing map,
critical points, potentials.
\end{abstract}

Existence of a Riemannian metric of constant positive curvature
in a given conformal class on a punctured sphere with conic
singularities at the punctures and prescribed angles at the
singularities is an important problem, but at present the solution
in this generality seems to be out of reach, see for example the
surveys in the introductions of \cite{EGT,MP}.

Mondelo and Panov \cite{MP} proposed a reduced problem:
to describe possible angles at the singularities
for such metrics.
(The conformal class is not prescribed.)
They solved this problem for generic angles.
In this paper their solution
is completed by the study of the remaining case
which was excluded in \cite{MP}.

The results in \cite{MP} are the following.
Let
\begin{equation}\label{A}
{\bm{\alpha}}=\{\alpha_1,\ldots,\alpha_n\},\quad
\alpha_j>0,\quad\alpha_j\neq 1
\end{equation}
 be the set of angles.
{\em We measure all angles in turns: $1$ turn is $2\pi$ radians.}
Strictly speaking, $\bm{\alpha}$ his is an unordered {\em multiset};
some elements can be repeated. Alternatively, ${\bm{\alpha}}$ is an
element of the $n$-th symmetric power of $R_{>0}$.
Whenever it is convenient, we list non-integer $\alpha_j$'s first,
followed by integer $\alpha_j$'s.

The restrictions are: the Gauss--Bonnet theorem,
\begin{equation}\label{GB}
\sum_{j=1}^n(\alpha_j-1)+2>0,
\end{equation}
and
\def\Z{\mathbf{Z}}
\begin{equation}\label{H}
d_1(\Z^n_o,{\bm{\alpha}}-{\mathbf{1}})\geq1,
\end{equation}
where ${\mathbf1}=(1,\ldots,1)$, and $\Z^n_o$ is the subset of the
integer lattice consisting of vectors with odd sums of coordinates,
and $d_1$ is the $\ell_1$ distance.

Mondello and Panov proved that these two conditions are always necessary,
and if one replaces (\ref{H}) by the strict inequality,
they also become sufficient.
So to complete their description, it remains to
investigate the case of equality in (\ref{H}):
\begin{equation}\label{H'}
d_1(\Z^n_o,\bm{\alpha-1})=1.
\end{equation}
Moreover, they proved that every metric satisfying (\ref{H'}) is
{\em co-axial}, which means that the monodromy group of the developing map
is a subgroup
of the unit circle. This gives a motivation for study of metrics with
co-axial monodromy. Here and in what follows, $S$
is the Riemann sphere (compact simply-connected Riemann surface).
The sphere with the standard spherical metric will be denoted by $\bC$.
In general, co-axial monodromy does not
imply~(\ref{H'}).

We say that a multiset (\ref{A}) is
{\em admissible} if there exists a metric of constant
curvature $1$ on $S\backslash\{ n\;\mathrm{points}\}$ with
conic singularities at these $n$ points with angles $\alpha_j$,
and the developing map of this metric has co-axial monodromy.
Our main result is
\vspace{.1in}

\noindent
{\bf Theorem 1.} {\em For a multiset (\ref{A}) assume that
$\alpha_{m+1},
\ldots,\alpha_n$ are integers while $\alpha_1,\ldots,\alpha_m$ are
not integers. For $\bm{\alpha}$ to be admissible it is necessary
that there exist a choice of signs $\epsilon_j\in\{\pm1\}$
and a non-negative integer $k'$ such that  
\begin{equation}\label{alt}
\sum_{j=1}^m\epsilon_j\alpha_j=k',
\end{equation}
and the number
\begin{equation}\label{bet}
k'':=\sum_{j=m+1}^n\alpha_j-n-k'+2\quad \mbox{is non-negative and even}.
\end{equation}
If the coordinates of the vector
\begin{equation}\label{c1}
\bm{c}:=(\alpha_1,\ldots,\alpha_m,
\underbrace{1,\ldots,1}_{ k'+k''\; \text{times}})
\end{equation}
are incommensurable, then (\ref{alt}) and (\ref{bet}) are also sufficient.

If $\bm{c}=\eta \bm{b}$, where coordinates of $\bm{b}$ are 
integers whose greatest common factor is $1$,
then there is an additional necessary condition
\begin{equation}\label{21}
2\max_{m+1\leq j\leq n}\alpha_j\leq \sum_{j=1}^q|b_j|.
\end{equation}

Conditions (\ref{alt}), (\ref{bet}), (\ref{21}) are sufficient
for $\bm{\alpha}$ to be admissible.}
\vspace{.1in}

As a corollary we mention that unless $n=2$, a co-axial metric
must have some integer angles whose sum is at least $n+k'-2$
and has the same parity as $n+k'$, where $k'$ is a number that
satisfies (\ref{alt}) with some choice of signs.
\vspace{.1in}

We briefly recall the definition of the
developing map of a metric of constant
positive curvature. Start with a small region
in $S\backslash\{{\mathrm{singularities}}\}$. It is well-known
that there is an isometry
from this region to a region in the standard sphere $\bC$. This isometry is conformal and thus analytic, and it admits an analytic continuation
along every curve which does not pass through the singularities.
We obtain a multi-valued function
$f:S\backslash\{{\mathrm{singularities}}\}\to\bC$ (or a genuine function on
the universal covering) which is called the developing map.
Conic nature
of the singularities means that $f(z)=f(a)+(c+o(1))z^\alpha$
near a singularity $a$ with angle the $\alpha$,
or $f(z)=(c+o(1))z^{-\alpha}$
if $f(a)=\infty$, where $z$ is a local conformal coordinate 
which equals $0$ at $a$, and $\alpha$ is the angle,
$f(a)$ means the
radial limit when $z\to 0$, and $c\neq 0$.
The result $f_\gamma$ of an analytic continuation of
$f$ along a closed path $\gamma$ not passing through the singularities
is related to the original germ of $f$
by $f_\gamma=\phi\circ\gamma$, where $\phi$ is an isometry of $\bC$,
so we obtain a representation of the fundamental group
of $S\backslash\{{\mathrm{singularities}}\}$ in the group of 
linear-fractional transformations.
The image of this representation is called
the monodromy group, and the developing map and the metric are called
co-axial if this monodromy group is a subgroup of the unit circle.
See \cite{E,EGT,MP}.

Suppose now that $f:S\to\bC$ is the developing map of
a metric with co-axial monodromy.
Monodromy group consists of transformations $w\mapsto\lambda w,\;
|\lambda|=1$. Then the meromorphic
$1$-form
$df/f$ is well defined on the sphere,
so 
\begin{equation}\label{R}
R=f'/f
\end{equation}
is a rational function. We assume without loss of generality
that $\infty\in S$ is not a pole of $R(z)dz$, so that
$R$ has a zero of order two at $\infty$.
{}From the local considerations we see that every pole of $R$ is simple,
the residue $\beta$ is real,
and every pole of $R$ is a conic singularity of the metric
with the angle $\alpha_j=|\beta|$,
unless $\beta=\pm1$. In the last case the pole of $R$
is a non-singular point
of the metric.
Moreover, a zero of $R$ in $\C$ of multiplicity $r$ is a conic
singularity
with the angle $\alpha_j=r+1$.
As $R$ has a double zero at infinity, the number of zeros (counting
multiplicity) in $\C$
is $q-2$, where $q$ is the number of poles.

Thus
\begin{equation}\label{R1}
R(z)=\sum_{j=1}^{m}\frac{\epsilon_j\alpha_j}{z-a_j}+\sum_{j=1}^{k}\frac{\delta_j}{z-b_j},
\end{equation}
where $\epsilon_j\in\{\pm1\},\;\delta_j\in\{\pm1\},$ and $R$ has zeros
in $\C$ whose sum of multiplicities is $q-2$, $q=m+k$,
and these multiplicities
are $\alpha_j-1$ for $m+1\leq j\leq n$. 
The developing map itself is thus given by
\begin{equation}\label{dev}
f(z)=\prod_{j=1}^q(z-z_j)^{\beta_j},
\end{equation}
where $\beta_j=\epsilon_j\alpha_j$ for $1\leq j\leq m$ and
$\beta_j=\delta_{j-m}$ for $m+1\leq j\leq q$.
The condition
\begin{equation}\label{cond}
\sum_{j=1}^{m}\epsilon_j\alpha_j+\sum_{j=1}^{k}\delta_j=0
\end{equation}
holds by the residue theorem.
As the number of zeros of $R$ in $\bC$ must be $m+k-2$,
and each zero is a singularity of the metric,
we obtain
\begin{equation}\label{a}
\sum_{j=m+1}^n(\alpha_j-1)=m+k-2.
\end{equation}
So a necessary condition for $\bm{\alpha}$ to be admissible is
\vspace{.1in}

\noindent
{\bf Condition 1.} {\em There exists a partition $\bm{\alpha}=A\cup B$
into two sub-multisets $A=\{\alpha_1,\ldots,\alpha_m\}$ and
$B=\{\alpha_{m+1},\ldots,\alpha_n\}$ so that:

All elements of $B$ are integers, and

There exist an integer $k$, and a choice of signs
$\epsilon_j\in\{\pm1\},\; 1\leq j\leq m$ and $\delta_j\in\{\pm1\},\; 1\leq j\leq k$,
such that (\ref{cond}) and (\ref{a}) hold.}
\vspace{.1in}

\noindent
{\bf Proposition 1.} {\em Condition 1 implies (\ref{alt}) and (\ref{bet}).}
\vspace{.1in}

{\em Proof.} Condition 1 coincides with (\ref{alt}) and (\ref{bet}) when
$A$ contains no integers. In this case we have $k=k'+k''$.

If for some $\bm{\alpha}$ a partition $A,B$
and numbers $k,\epsilon_j,\delta_j$ as required by Condition 1 exist,
then there exists another
partition $A',B'$ with the same properties and with the additional
property that $A'$ contains no integers.

Indeed, suppose that $A=\{\alpha_1,\ldots,\alpha_m\}$ and $\alpha_m$
is an integer. Then define $A'=\{\alpha_1,\ldots,\alpha_{m-1}\}$
and $B'=\{\alpha_m,\ldots,\alpha_n\}$. To restore (\ref{cond}) we must
add $\alpha_m$ of $\delta_j$'s equal to $-\epsilon_m$; this increases $k$
to $k^*=k+\alpha_m$ and decreases $m$ to $m^*=m-1$,
so the total increases
in the right and left hand sides of (\ref{a}) are equal,
so this
condition (\ref{a}) is satisfied for the new partition $A',B'$.
We repeat this procedure until all integer angles are removed
from $A$. In the special case when all angles are integers,
$A$ will be empty.
This proves the proposition and necessity of conditions
(\ref{alt}) and (\ref{bet}) in Theorem 1.
\vspace{.1in}

We will call such partitions where $A$ consists of all non-integer angles
of $\bm{\alpha}$ {\em reduced}. 
In a reduced partition of $\bm{\alpha}$,
$m$ is the number of non-integer angles, and the only reason why
a reduced partition may be non-unique is that different choices
of signs $\epsilon_j$ in (\ref{R1}) may be possible.

When the number of non-integer angles $m\leq 3$,
conditions equivalent to (\ref{alt}) and (\ref{bet})
were obtained in \cite{E}, \cite[Thm. 4.1]{EGT},
\cite{EGT2}, 
and for $m\leq 3$ they are also sufficient. 

Formula (\ref{R1}) for a reduced partition can be written as
$$\frac{f'}{f}=\sum_{j=1}^m\frac{\epsilon_j\alpha_j}{z-a_j}-
\sum_{j=1}^{k'}\frac{1}{z-b_j}+
\sum_{j=k'+1}^{k'+k''}\frac{(-1)^j}{z-b_j},$$
and the residue theorem combined with (\ref{alt}) shows that $k''$
must be even,
as stated in (\ref{bet}).

We will see that for co-axial metrics, (\ref{alt}) and (\ref{bet})
are also sufficient in the generic situation, when the coordinates
of the vector $\bm{c}$ in (\ref{c1}) 
are incommensurable. When coordinates of $\bm{c}$
are commensurable, there are additional restrictions. 

For a given multiset
${\bm{\alpha}}=\{\alpha_1,\ldots,\alpha_n\}$ satisfying
Condition 1 
we call a quadruple $(A,B,k,\{\epsilon_j\},\{\delta_j\})$ of
parameters in (\ref{R1}) an
{\em arrangement} for ${\bm{\alpha}}$. For an admissible
$\bm{\alpha}$, different arrangements may give
different metrics. We do not require that
all $\alpha_j\in A$ are non-integers. If they are non-integers,
the arrangement is called {\em reduced}.
If ${\bm{\alpha}}$ is admissible, there exist finitely many
arrangements, at least one of them
is reduced. This reduced arrangement
may be non-unique: various choices of $\epsilon_j$ in (\ref{alt}) are
sometimes possible. A priori, we have to deal
with non-reduced arrangements because we have further
conditions besides conditions (\ref{alt}) and (\ref{bet});
it is possible that some of the arrangements
satisfy them, others do not.

The geometric meaning of a reduced arrangement is the following.
The developing map $f$ as in (\ref{dev}) is a multi-valued function, but
the preimage $f^{-1}(\{0,\infty\})$ is well defined (as radial limits).
Metrics corresponding
to reduced arrangements are exactly those for which the
developing map does not take the values $0,\infty$ at the singularities
with integer angles. 

{}From now on, we assume that Condition 1 is satisfied,
an arrangement (perhaps not reduced)
is fixed, and the logarithmic derivative of the developing map
is written as in (\ref{R1}).
We denote for simplicity
\begin{equation}\label{c}
\{\epsilon_1\alpha_1,\ldots,\epsilon_m\alpha_m,\delta_1,\ldots,\delta_k\}=
\{ c_1,\ldots,
c_q\},
\end{equation}
where $q=m+k$. In the case of reduced arrangement $k=k'+k''$.
The question is what multiplicities of zeros of $R$ are possible
for a given 
vector of residues $c$:
\vspace{.1in}

\noindent
{\bf Question 1.} {\em Suppose that real non-zero numbers
$\{ c_1,\ldots,c_q\}$ are given, and 
\begin{equation}\label{sum}
\sum_{j=1}^q c_j=0.
\end{equation}
Which partitions
$\{ \ell_1,\ldots,\ell_k\}$ of $q-2$ can be realized as multiplicities
of zeros in $\C$ of the function
\begin{equation}\label{g}
g(z)=\sum_{j=1}^q\frac{c_j}{z-z_j},
\end{equation}
where $z_j$ are pairwise distinct complex numbers?}
\vspace{.1in}

Notice that zeros of $g$ are critical points of the potential
\begin{equation}\label{potential}
u(z)=\sum_{j=1}^qc_j\log|z-z_j|,
\end{equation}
and Question 1 seems to be of independent interest.

The trivial but important property is the following:
\vspace{.1in}

{\em If all $c_j$ are multiplied by a constant, the multiplicities
of zeros of $R$ do not change.}
\vspace{.1in}

So we introduce the real projective space $\RP^{q-2}$ which consists
of non-zero $q$-tuples $(c_1,\ldots,c_q)$ satisfying (\ref{sum}), modulo
proportionality. A point $\bm{c}\in\RP^{q-2}$ is called {\em rational}
if its equivalence class contains a $q$-tuple with all $c_j$ rational.
Let $Z$ be the union of coordinate hyperplanes
$Z_j=\{ c:c_j=0\},\; 1\leq j\leq q$.
Let $P$ be a partition of $q-2$. A point $\bm{c}\in\RP^{q-2}$
is called {\em $P$-admissible},
if there exist pairwise distinct $z_j\in\C$
such that the function $g$ in (\ref{g}) has zeros of multiplicities $P$.
A point $\bm{c}\in\RP^{q-2}\backslash Z$
which is not $P$-admissible is called {\em $P$-exceptional}.
\vspace{.1in}

\noindent
{\bf Theorem 2.} {\em Let $P=\{\ell_1,\ldots,\ell_s\}$
be a partition of $q-2$.
Every irrational point $\bm{c}\in\RP^{q-2}\backslash Z$ is $P$-admissible.
A rational point if $P$-admissible if and only if
\begin{equation}\label{proh}
2(1+\max_{1\leq j\leq s}\ell_j)\leq \sum_{j=1}^q|b_j|,
\end{equation}
where
\begin{equation}
\label{b}
\{ b_1,\ldots, b_q\}=\{\eta c_1,\ldots,\eta c_q\},\quad\eta\neq 0,
\end{equation}
is a vector with mutually prime integer coordinates proportional
to $\bm{c}$.}
\vspace{.1in}

When all residues in (\ref{g}) are mutually prime integers 
the developing map $f=\exp\int g$
is a rational function,
for which positive residues are multiplicities of zeros, negative
residues are multiplicities of poles, and $k_j:=\ell_j+1,
\;1\leq j\leq s$ are the multiplicities of $f$ at other critical
points, different from zeros and poles. So our Question 1
with mutually prime integer residues is
a special case 
of the Hurwitz problem \cite{B,M}:
\vspace{.1in}

\noindent
{\bf Question 2.} {\em Given two partitions
$\{ n_1,\ldots,n_r\}$ and $\{ m_1,\ldots,m_t\}$ of the same number $d>1$
and a multiset of integers $\{ k_1,\ldots,k_s\},\; k_j\geq 2,$
such that 
\begin{equation}\label{RH}
\sum_{j=1}^r(n_j-1)+\sum_{j=1}^t(m_j-1)+\sum_{j=1}^s(k_j-1)=2d-2,
\end{equation}
does there exist a rational function $f$ of degree $d$
with zeros of multiplicities $m_j$ and poles of multiplicities
$n_j$ and other critical points were multiplicities of $f$ are
$k_j$?}
\vspace{.1in}

Here $m_j$ are positive residues in (\ref{g}) and $n_j$
are negative residues in (\ref{g}). 

By a simple perturbation argument, it is sufficient to consider the
case when the values of $f$ at these ``other critical points''
are all distinct: critical points with the same critical
value other than $0,\infty$ can be perturbed so that all these
critical points will have different critical values, and the
multiplicities of zeros and poles are not affected.

The answer to Question 2 was
recently obtained by Song and Xu \cite{SX}. 
The necessary and sufficient condition of existence of $f$ is
\begin{equation}\label{cond2}
k_j\leq d=\frac{1}{2}\sum_{i=1}^q |b_i|,\quad 1\leq j\leq s.
\end{equation}
The necessity of this condition is evident because
the right hand side is the degree $d$ of $f$. Sufficiency was proved
by Song and Yu who generalized the result of
Boccara \cite{B} for $s=1$.
This solves Question 2 and proves Theorem 2 for the case of
rational vector $\bm{c}$.

We state a trivial but important 
\vspace{.1in}

\noindent
{\bf Remark 1.} {\em For each $P$ there are finitely many
integer vectors $(b_1,\ldots,b_q)$ which do not satisfy (\ref{cond2}),
thus there are finitely many rational points in $\RP^{q-2}\backslash Z$
which are $P$-exceptional.}
\vspace{.1in}

Theorem 2 gives an 
algorithm which determines
whether a given multiset $\bm{\alpha}$ is admissible.
The algorithm works as follows. Starting with a multiset
${\bm{\alpha}}=\{\alpha_1,\ldots,\alpha_n\}$ we check conditions
(\ref{alt}) and (\ref{bet}).
If they are not satisfied, then $\bm{\alpha}$ is not admissible.
If these conditions are satisfied, we consider all arrangements
for $\bm{\alpha}$ and vectors $\bm{c}$ corresponding to them as in
(\ref{c}). If one of these vectors is irrational,
then $\bm{\alpha}$ is admissible. If all are rational,
we construct integer vectors $\bm{b}$ as in (\ref{b}). If one
of these vectors $\bm{b}$ satisfies (\ref{proh})
with $\ell_j=\alpha_{j+m}-1$,
 then $\bm{\alpha}$
is admissible, if none, then not.  
\vspace{.1in}

Most of Theorem 1 is a corollary of Theorem 2, except the statement
that it is enough to check condition (\ref{21}) only for one
reduced arrangement. This will be addressed in the formal proof of
Theorem 1 in the end of the paper.
\vspace{.1in}

\noindent
{\bf Example 1.} For arbitrary non-integer
$\beta>0$ the multiset $\{\beta,\beta,\beta,\beta,3\}$ is not admissible.
Conditions (\ref{alt}) and (\ref{bet}) are satisfied.
The only reduced arrangement is $A=\{\beta,\beta,\beta,\beta\},\;
B=\{3\},\; k=0.$ So $q=4$, $b=(1,1,-1,-1)$,
and condition (\ref{21}) is
violated.
\vspace{.1in}

\noindent
{\bf Example 2.} For arbitrary $\beta>0$ the
multiset $\{\beta,\beta,2\beta,2\beta,3\}$ is admissible.
Conditions (\ref{alt}) and (\ref{bet}) are satisfied.
Take the arrangement
$A=\{\beta,\beta,2\beta,2\beta\},\; B=\{3\},\; k=0.$
We have $q=4$, $b=(1,-1,2,-2)$. 
Inequality (\ref{21}) is satisfied.
Let us write a developing map explicitly for this case:
$$f(z)=\left(\frac{(z-1)^2(2z+1)}{(z+1)^2(2z-1)}\right)^\beta
=h^\beta(z)$$
The corresponding metric has angles $\beta$ at $\pm 1/2$
and angles $2\beta$ at $\pm1$. In addition to this, there is
angle $3$ at $0$, because $h$ has a triple point at $0$
with critical value $-1$.
\vspace{.1in}

{\em Sketch of the proof of Theorem 2.} First we notice that
the problem of constructing a surface of constant positive curvature,
with co-axial monodromy and
with prescribed angles at conic singularities is equivalent to
a similar construction problem for a surface with a flat metric.
Trying to construct this flat surface by gluing cylinders, we discover
the general nature of obstructions: the given angles
must satisfy some systems of inequalities. These inequalities
are too complicated to write explicitly,
but we determine their general nature:
they are inequalities between some {\em linear forms} in the residues $c_j$
with {\em integer} coefficients. 
Therefore, for each
partition $P$, the set of $P$-exceptional
points is a rational polyhedron
in the space $\RP^{q-2}$. If this polyhedron
consists of infinitely many points, then it must also contain infinitely
many rational points. But we know from Theorem 2
that the number of exceptional
rational points is finite for given $P$, see Remark 1. 
Therefore the polyhedron of $P$-exceptional points consists
of finitely many points and {\em thus} all exceptional
points must be rational. 
\vspace{.1in}

{\em Proof of Theorem 2.}
\vspace{.1in}

\noindent
1. {\em From spherical to flat and back.}
\vspace{.1in}

Let $f:S\to\bC$ be the developing map as in (\ref{dev}). 
Let $\Omega=S\backslash\{ z_1,\ldots,z_q\}$. Then we have the restricted
map $f^*:\Omega\to\C^*$.

We equip $\C^*$ with the flat metric whose length element is
$|dz/z|$. This metric makes $\C^*$ into an open cylinder infinite
in two directions 
whose {\em girth} (the length of the shortest non-trivial geodesic,
a. k. a. the systole)
is $2\pi$.
We pull back this flat metric to $\Omega$ via $f^*$ and obtain a flat
surface
which is conformally equivalent to a sphere with $q$ punctures, and some
neighborhoods
of the punctures are semi-infinite cylinders of girths $2\pi|\beta_j|$.
We call this surface $(\Omega,\rho)$, where $\rho$ is the flat metric.
The developing map $f^*$ of $(\Omega,\rho)$ has two special features:
it maps $\Omega$ to $\C^*$ (rather then $\C$) and it tends to
$0$ or to $\infty$ at the punctures in the sense of radial limits.

Conversely, suppose that $\Omega$ is a Riemann surface
conformally equivalent
to a punctured sphere, equipped with a flat Riemannian metric $\rho$
 such that
some neighborhoods of the punctures are semi-infinite cylinders of
girth $2\pi\mu_j$. Moreover, suppose that the developing map
$h$ maps $\Omega$ to $\C^*$ and tends at each puncture either to $0$
or to $\infty$.
By filling the punctures,
we can extend $h$ to a (multivalued) map $f:S\to\bC$
and pull back the spherical metric to $S$. The resulting surface
has constant curvature $1$,
and in addition to conic singularities
in $\Omega$ has conic singularities at
the punctures $S\backslash\Omega$. The angles at
these additional singularities $z_j$
are $\mu_j$. 
\vspace{.1in}

\noindent
{\em 2. From flat surface to a system of linear inequalities.}
\vspace{.1in}

Now we study this auxiliary flat
surface $(\Omega,\rho)$ and its developing map $h$.
The level sets
$$L_t=\{ z\in\Omega:\log|h(z)|=t\},\quad -\infty<t<\infty,$$
make a foliation of $\Omega$.
This means that $\Omega$ is a disjoint union of
{\em leaves} and finitely many {\em critical points} of $\log|h|$.
Leafs are the curves on which $|h(z)|$
is constant; these curves are either simple closed curves (ordinary
leaves)
or simple open curves with both ends at singular points (singular leaves).
Foliations are considered here as topological objects: up to
homeomorphisms which respect leaves.

Notice that unlike the developing map $h$, the function $u=\log|h|$ is
a well-defined (single-valued) harmonic function.
Level sets $L_t$ which contain singular points are called
{\em critical level sets}. A non-critical level set consists
of finitely many ordinary leaves, while a critical level set
may contain both ordinary and singular leaves and some critical points.

The region $\Omega$ is a disjoint union of
open {\em foliated cylinders} and critical level sets.
A model foliated cylinder is obtained by taking a rectangle in the plane
foliated into horizontal segments and identifying its vertical sides
in the natural way. An open foliated
surface homeomorphic to such a cylinder,
by a homeomorphism respecting the foliation is called a foliated cylinder.

Every singular point in $\Omega$ is a saddle point of $u$ and
it has an index: a positive integer $k$
such that the singular leaves in a neighborhood
of this point look like
the $2(k+1)$ intervals of the set $\{ z:|z|<1,\,\Rea z^{k+1}=0\}$
meeting at $0$.
This is because our  function $u$ is harmonic.

Our foliation 
has an additional structure: there are two functions
on the set of leaves: one is the {\em height} $t$, another is the
length of a leaf with respect to the 
intrinsic metric $\rho$. For a leaf $\gamma\subset L_t$ the height is $t$.
The length of a leaf $\gamma$
is a positive number which can be computed
by the formula
$$|\gamma|=
\int_\gamma \left|\frac{\partial u}{\partial n}\right||dz|,\quad
u(z)=\log|h(z)|,$$
where $n$ is the unit normal to $\gamma$. The same
formula defines the length of any arc of a leaf.

Suppose that an interval $(t',t'')$ contains no critical values of $u$.
For $t\in (t',t'')$, let $\gamma_t\subset L_t$ be a leaf which depends
continuously on $t$. (Convergence of leaves which is used here
is uniform, using some parametrization).

Then {\em the length $|\gamma_t|$ does not depend on $t$.}

This follows from Green's formula applied to $u$ in the ring
between $\gamma_{t_1}$ and $\gamma_{t_2}$ where $t_1,t_2$ are
any numbers between $t'$ and $t''$.

When $t$ passes through a critical value, some leaves break into
singular leaves and then these singular leaves re-assemble into new
ordinary leaves.

More precisely, let $(t',t'')$ be as above, and suppose
that $t'$ is a singular value. Choose a leaf $\gamma_t\in L_t$
which depends continuously on $t$ for $t\in(t',t'')$.
Then as $t\to t'+$ some parametrization of $\gamma_t$
converges uniformly to a closed curve, which can be an ordinary leaf,
or a finite union of singular leaves $\gamma^j\subset L_{t'}$ and singular points.
Moreover, we have
\begin{equation}\label{lin1}
|\gamma_t|=\sum_j|\gamma^j|,
\end{equation}
where the summation is over all those leaves which form the limit
of $\gamma_t$, and all summands
in the right hand side are strictly positive.

Relations (\ref{lin1}) form a system of linear equations which the
lengths of leaves of a given topological foliation must satisfy,
assuming that the lengths of ordinary leaves do not change with height.
\vspace{.1in}

\noindent
{\em 3.
From foliations with height and length back to flat surfaces.}
\vspace{.1in}

Suppose now that $\Omega$ is a topological punctured sphere with
a topological foliation whose leaves are level sets of
some smooth function $v:\Omega\to\R$ with finitely many critical points,
and $v(z)\to\pm\infty$
when $z$ tends to a puncture, and
\vspace{.1in}

\noindent
a) In a neighborhood of each critical point $v$
is topologically equivalent to a harmonic
function.
\vspace{.1in}

Suppose further that 
a strictly positive function $\phi$ on the set of leaves is given
which has the formal properties of the length function, namely:
\vspace{.1in}

\noindent
b) If $\gamma_t\subset L_t$ is a family of
ordinary leaves continuously depending on $t\in(t',t'')$
on an interval containing no critical values, then
$\phi(\gamma_t)$ is constant on $(t',t'')$.
\vspace{.1in}

\noindent
c) If $t'$ is a singular value, and $\gamma_t$ is the same as in b),
and $\gamma_t$ tends to the union of singular leaves $\cup_j\gamma^j$
as $t\to t'$, we have 
\begin{equation}\label{lin2}
\phi(\gamma_t)=\sum_j\phi(\gamma^j).
\end{equation}

We claim that whenever such a foliation and functions $t$ and 
$\phi$ on the leaves are given, one can
introduce a flat metric on $\Omega$ whose developing map
$h$ has the property
that the level sets of $v=\log|h|$ define our given foliation. 
and the function $\phi$ is the length of the leaves of this foliation.
The flat metric defines on $\Omega$ the conformal structure of
a punctured sphere.

To prove the claim, we consider the partition of $\Omega$ into foliated
cylinders
$C_j$ and critical level sets as described in part 2 of the proof.
Each cylinder is mapped by $v$
into a
maximal interval $(t',t'')$ free of
critical values of $v$.  In the {\em trivial case} when there are no
critical points at all, we have $(t',t'')=(-\infty,\infty)$.
In all other cases
there are two such semi-infinite intervals and finitely many finite
intervals.

Each foliated cylinder $C_j$ is homeomorphic to
the product $\gamma_j\times (t',t'')$, where $\gamma_j$ is an ordinary
leaf in $L_t$ for some $t\in(t',t'')$.
We pull back to $C_j$ the standard Euclidean metric from this product,
so that $|\gamma_j|=\phi(\gamma_j)$.
This defines the flat metric $\rho$ on the cylinders of the foliation.
Some of them are of finite height, others semi-infinite,
except the trivial case when there is only one doubly-infinite cylinder.

Let $\overline{C_j}$ be completions of the $C_j$ with respect to their
metrics. The boundary circles of $C_j$ correspond to some leaves of
the foliation on the singular level sets, and some finite sets of points
on each boundary circle must be glued together into singular points.
So we break every boundary circle into arcs which will correspond
to the singular leaves. The lengths of these arcs
are determined by our function $\phi$, and this is where
relation (\ref{lin2}) is used. Then we glue together our cylinders
along these arcs
respecting the length. To perform this gluing we use the theorem of
Aleksandrov and Zalgaller, see, for example \cite[Thm. 8.3.2]{Resh},
about gluing two surfaces along a geodesic arc.
It guarantees 
that we obtain a ``surface of bounded curvature'' in the sense
of Aleksandrov, with a flat metric and finitely many conic
singular points.
(This is a surface of special kind which is called
a polyhedral surface in \cite{Resh}).
The total angle at a singularity is equal to
one half of the number of boundary points of
cylinders which are glued together at this point.

That the resulting surface is connected and of genus $0$ is guaranteed
by the topology of the foliation. For a given foliation, the only
condition for the possibility of this gluing is the linear relations (\ref{lin2})
between the lengths of the leaves. 
\vspace{.1in}

\noindent
{\em 4. Conclusion of the proof of Theorem 2.}
\vspace{.1in}

Suppose that the vector $c$ has $p$ positive and $r$ negative coordinates,
$p+r=q$. Take a $q$-punctured sphere $\Omega$, and construct
a function $v:\Omega\to\R$ which
tends to $-\infty$ at $p$ punctures of $\Omega$ and to $+\infty$
at the remaining $r$ punctures. Moreover, we require that
all critical points of $v$ in $\Omega$ are saddle points of the
topological types
which are possible for harmonic functions, and the multiplicities
of these critical points are the parts of the partition $P$.
\vspace{.1in}

\noindent
{\bf Lemma.} {\em For every $p$, $r$ and $P$ there exists a function $v$
with these properties.}
\vspace{.1in}

Postponing the proof of the Lemma, we complete the proof of Theorem~2.

Consider {\em all} foliations defined by functions $v$ satisfying our
conditions with fixed $p,r,P$.
To assign a length function $\phi$ consistent
with a foliation, we have to solve the system of equations (\ref{lin2})
which is determined by the foliation. In this system, the given numbers
are the girths of the semi-infinite cylinders (these are our $|c_j|$),
and the unknown variables are the girths of all finite height
cylinders and the lengths of the singular leaves.

In addition to (\ref{lin2}),
we have the restriction that all girths and lengths
must be strictly positive.

If this linear system has a strictly positive solution, we can construct
our metric by performing steps described in parts 3 and 1.
If not, a metric
corresponding to this particular foliation does not exist.

We give an illustrating example. Suppose we want to construct
a function $g$ as in (\ref{g}) with two positive residues $a,b$,
two negative residues $-c,-d$ and a single critical point where
the local degree of $h$ is $3$. Then the critical level set must have the
form as in Fig.~1, where the regions represent semi-infinite cylinders,
the dots labeled $a,b$ are the poles of $h$, and the dots labeled
$c,d$ are zeros of $h$.
The girths of the four cylinders are
$2\pi a,2\pi b,2\pi c,2\pi d$ and the length of the three singular leaves are
$2\pi b$, $2\pi d$, and  $2\pi x$ (see Fig.~1).
Non-singular leaves of the foliation
are not shown
in the picture, each of them is a Jordan curve that surrounds
the puncture in its cylinder.
Equations (\ref{lin2}) for this case are 
$$a=x+d,\quad c=x+b,$$
which are consistent if and only if $a-d+b-c=0$.
Now $x$ must be strictly positive, so we obtain necessary
and sufficient conditions
of existence of such $g$: $a>d$ and $c>b$, in other words, the
positive residues must be unequal and negative residues must be unequal.
This explains examples 1,2 above.
\begin{figure}
\centering
\includegraphics[width=1.5in]{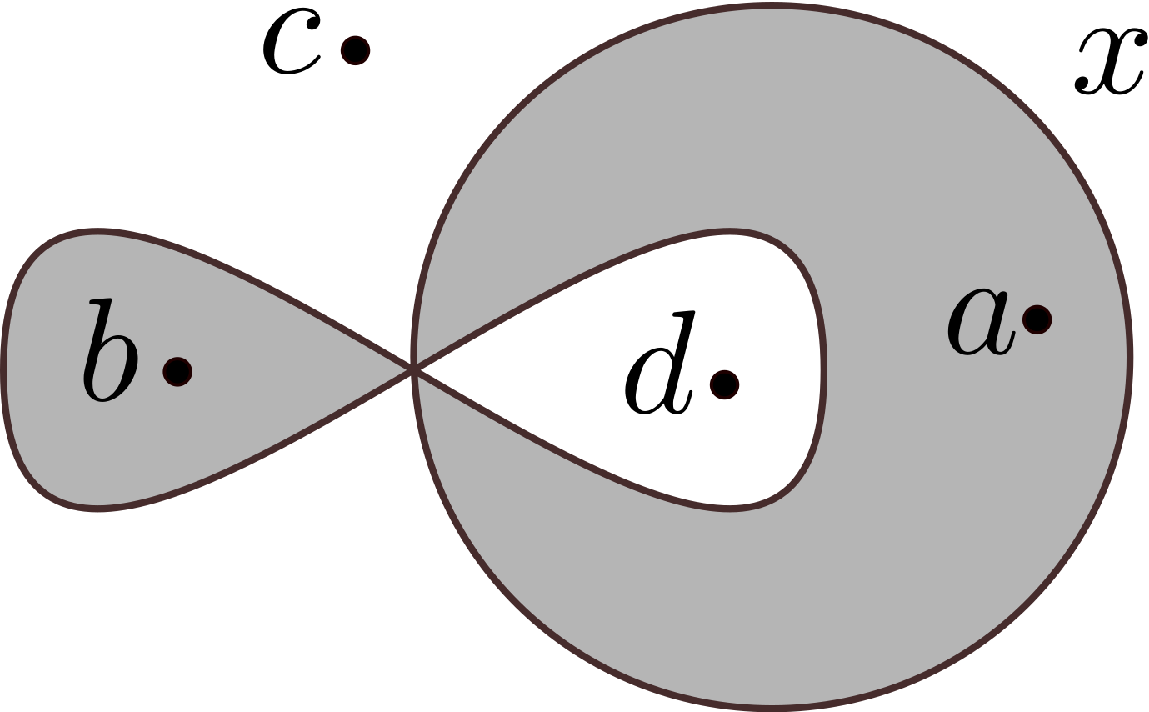}
\caption{A critical level set.}
\end{figure}

In any case, the condition that a vector $c$ is $P$-admissible is
stated in terms of linear equations and linear inequalities with
integer coefficients and Boolean operations.
So $P$-exceptional vectors $c$ form a rational
polyhedron in $\RP^{q-2}$. If this polyhedron is infinite, then
it contains
infinitely many rational points \cite{D}, which is not the case:
we have seen that there are only finitely many $P$-exceptional rational
vectors for each given $P$  (see Remark 1). 
So the polyhedron is finite.
So it consists of only rational points.

This completes the proof of Theorem~2.
\vspace{.1in}

\noindent
{\em 5. Proof of Lemma 1.}
\vspace{.1in}

It is sufficient to prove the lemma for the special case when there is
only one critical point. Then it can be broken into pieces according
to partition $P$ by a perturbation, as shown in the first three
lines of figure 2. In lines 1-3, on the left hand side we have
a critical point of multiplicity $4$.
In lines 1 and 2 it is broken to two critical points of
multiplicity $2$, in line 3 it is broken into
one critical point of multiplicity $2$ and two critical points
of multiplicity $1$.
\begin{figure}
\centering
\includegraphics[width=4.5in]{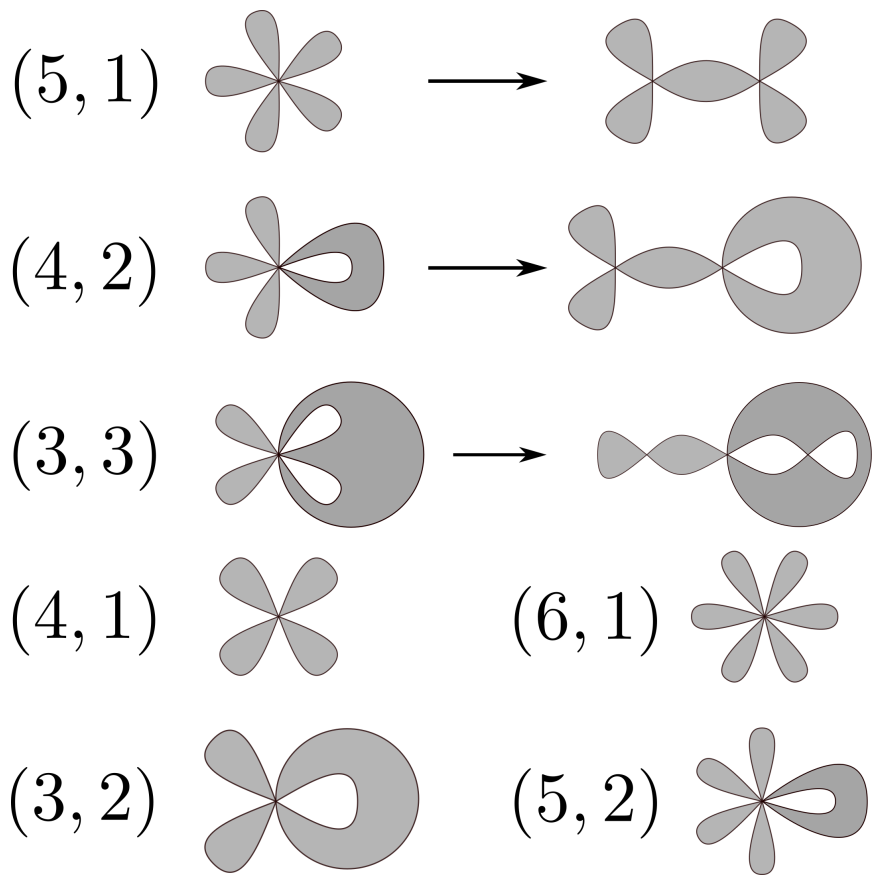}
\caption{
Black and white regions represent semi-infinite cylinders,
numbers on the left are $(p,r)$, and the arrows
arrows show breaking a high multiplicity critical point into
critical points of the lower multiplicity.}
\end{figure}
A foliation with one critical point is defined by its critical level set,
say $v(z)=0$,
and by
assigning a black or white color to the components of the complement
according to the sign of $v$. Instead of describing these foliations
in words we just present pictures of their critical level sets in
Fig.~2 . Each black or white region represents a semi-infinite
cylinder with one puncture inside. Numbers $(p,r)$ are written
on the left. We start with $(p,1)$, a flower with $p$ petals,
and then pass to $(p-1,2)$, $(p-2,3)$ etc, as shown in the picture.
\vspace{.1in}

\noindent
{\em Proof of Theorem 1.}
\vspace{.1in}

Necessity of conditions (\ref{alt}), (\ref{bet}) has been already
explained, and sufficiency
of (\ref{alt}), (\ref{bet}) and (\ref{21}) follows from Theorem 2. 
It remains to prove the necessity of (\ref{21}). It is necessary
that (\ref{21}) is satisfied for some arrangement. We have to prove
that it is enough to check it only for one reduced arrangement.

Suppose that $\bm{\alpha}$ is a multiset satisfying Condition 1.
We claim that if it is admissible, then there exists a metric
with these angles corresponding to some reduced arrangement.
Indeed if there are singular points with integer angles
for which the developing map takes the values $0$ or $\infty$,
then one can find another metric with co-axial monodromy
with the same angles for which the developing map does not take
the values $0$ or $\infty$ at the singular points with integer angles.
This follows from a general argument which permits to
``move around'' a singular point with integer angle.

Let $f:S\to\bC$ be the developing map of
a surface of curvature $1$ with conic singularities,
and suppose that $a\in S$
is a singular point with integer angle $\alpha$.
Let $r>0$ be smaller than
the distance from $a$ to other singularities, and such that
the closed intrinsic disk $D$
of radius $r$ centered at $a$ is homeomorphic to a closed disk
in the plane.

We will remove the interior of $D$ from $S$,
\def\int{\mathrm{int}\, }
and paste $S\backslash\int D$ with a new
surface $C$ homeomorphic to a closed disc in the plane,
equipped with a metric
of the same constant curvature, having one singularity in the interior
with the same angle $\alpha$. This can be so arranged that
the distance in $C$ from the singularity to $\partial C$ is any positive number
less than $r$, and the closest point to the singularity on $\partial C$
is any given point of $\partial C$.
So we have a continuous family of
deformations. Moreover, the resulting surface
$S'=C\cup (S\backslash\int D)$
is smooth, and has constant curvature except at the conic singularities
in $S\backslash D$ and in $C$.

Consider the disk $U=\{ z:|z|<R\},\; R=\tan(r/2)\}$ equipped
with the standard spherical metric $\rho$. (The spherical radius
of this disk is $r$.)
Let $C=\{ z:|z|\leq R^{1/\alpha}\}$ equipped with the metric
$\rho_1=f^*\rho$, where $f(z)=z^\alpha$.
Let $D=\{ z:|z|\leq R^{1/\alpha}\}$ equipped
with the metric $\rho_2=g^*\rho$, where
$$g(z)=R\frac{z^\alpha+aR}{1+\overline{a}z^\alpha},$$
where $|a|<1$. Consider the annulus $B\subset U$,
$B=\{ z:t<|z|<R\}$ where $t\in(R|a|,R)$,
and let denote $A_1=f^{-1}(B),\; A_2=g^{-1}(B)$. Then the metric
spaces $(A_1,\rho_1)\subset D$ and $A_2\subset C$ are isometric,
because they are both isometric to the covering of $B$ of degree $\alpha$.
So we can remove from our surface $S$ a disk isometric to $D$
and glue in $C$ instead. The parameter of deformation is $a$.
Notice that this deformation is {\em isomonodromic}, does not
change the monodromy group. 

This explains why the necessary condition (\ref{21}) in Theorem 1 is 
enough to verify for reduced arrangements: we can always perturb
a co-axial metric and obtain another co-axial metric
with the same angles and reduced arrangement.

Now we notice, that
if some reduced arrangement satisfies (\ref{21}) of Theorem 1,
then all other 
reduced arrangements for the same multiset of
angles will also satisfy (\ref{21}), because $q$ and $\sum|b_j|$
are the same
for all reduced arrangements. Indeed, $m$ and
$$q=m+k'+k''=\sum_{j=m+1}^n(\alpha_j-1)+2,$$
depend only on $\bm{\alpha}$, and $b_j$ depend only on
non-integer angles in $\bm{\alpha}$ 
and on $k=k'+k''$.
This proves necessity of condition (\ref{21}) and completes the proof
of Theorem~1. 
\vspace{.1in}

{\em Remark 2.} A similar deformation of a singularity with
non-integer angle
is impossible. Consider, for example a ``football'',
the sphere with
a metric of curvature $1$ and
two conic singularities. The singularities of such
surface must have equal
angles, and for each angle there is such a surface. But if the angle is
non-integer, then a football is unique, while with an integer
angle there is a $1$-parametric family of footballs \cite{Tro}.
What was used
in our argument is that the developing map is single-valued
in a neighborhood of a singularity with integer angle.
\vspace{.1in}

{\em Remark 3.} If $\bm{\alpha}$ is an admissible multiset, there exists
a metric of positive curvature with angles $\bm{\alpha}$.
But the conformal
class of this metric cannot be arbitrarily assigned. Take for example
$\bm{\alpha}=(\alpha_1,\ldots,\alpha_n)$ where $\alpha_n=n-2$, the
rest of the angles are not integers and
$$\sum_{j=1}^{n-1}\alpha_j=0.$$
Compare Example 2 above.
The developing map satisfies 
$$\frac{f'(z)}{f(z)}=\sum_{j=1}^{n-1}\frac{\pm\alpha_j}{z-z_j},$$
and the right hand side must have a zero of multiplicity $n-3$.
This imposes $n-4$ conditions on the poles $z_j$. Indeed,
we may assume without
loss of generality that $z_1=0,\; z_2=1$, so we obtain $n-4$ conditions
on $n-3$ variables $z_j$
which suggests that there is only a one-dimensional
family of conformal classes of such metrics.

Similar phenomenon may occur when all angles are non-integer.
Lin and Wang \cite{LW} studied a problem which is equivalent
to description of metrics of positive curvature on the sphere with
four singularities
with angles $(1/2,1/2,1/2,3/2)$. The conformal type of these metrics
depends on one complex parameter,
and it turns out that the moduli space of
quadruply punctured spheres is split into two parts,
each with non-empty interior,
such that for one part a metric with these angles exists and for 
the other part it does not. In all these examples the angles
are very special. The results in \cite{BMM,CL} suggest that perhaps
for generic angles satisfying (\ref{GB}) and (\ref{H}) a metric
of curvature $1$ exists in prescribed conformal class of the
punctured sphere. 
\vspace{.1in}

The author thanks Andrei Gabrielov, Michael Kapovich,
Dmitry Novikov, Carlo Petronio and Vitaly Tarasov
for helpful discussions.

{\em Department of Mathematics, Purdue University,

West Lafayette, IN 47907 USA

eremenko@math.purdue.edu}
\end{document}